# Generalized Bombieri – Lagarias' theorem and generalized Li's criterion

Sergey K. Sekatskii

*Laboratoire de Physique de la Matière Vivante, IPSB, Ecole Polytechnique Fédérale de Lausanne, BSP, CH 1015 Lausanne, Switzerland*
E-mail : Serguei.Sekatski@epfl.ch

We show that Li's criterion equivalent to the Riemann hypothesis, viz. the statement that the sums $k_n = \Sigma_\rho (1-(1-\frac{1}{\rho})^n)$ over Riemann xi-function zeroes as well as all derivatives $\lambda_n \equiv \frac{1}{(n-1)!}\frac{d^n}{dz^n}(z^{n-1}\ln(\xi(z)))|_{z=1}$, where $n=1, 2, 3...$, should be non-negative if and only if the Riemann hypothesis holds true, can be generalized and the non-negativity of certain derivatives of the Riemann xi-function estimated at an *arbitrary* real point $a$, except $a=1/2$, can be used as a criterion equivalent to the Riemann hypothesis. Namely, we demonstrate that the sums $k_{n,a} = \Sigma_\rho (1-\left(\frac{\rho-a}{\rho+a-1}\right)^n)$ for any real $a$ not equal to ½ as well as the derivatives $\frac{1}{(n-1)!}\frac{d^n}{dz^n}((z-a)^{n-1}\ln(\xi(z)))|_{z=1-a}$, $n=1, 2, 3...$ for any $a<1/2$ are non-negative if and only if the Riemann hypothesis holds true (correspondingly, the same derivatives when $a>1/2$ should be non-positive for this). Similarly to Li's one, the theorem of Bombieri and Lagarias applied to certain multisets of complex numbers, is also generalized along the same lines. The results presented here are now available as a paper in *Ukrainian Math. J.*, **64**, 371-383, 2014.

1. INTRODUCTION

In 1997, Li has established the following criterion equivalent to the Riemann hypothesis concerning non-trivial zeroes of the Riemann $\zeta$-function (see e.g. [1] for standard definitions and discussion of the general properties of this function) and now bearing his name (Li's criterion) [2]:

**Li's criterion.** *Riemann hypothesis is equivalent to the non-negativity of the following numbers*

$$\lambda_n \equiv \frac{1}{(n-1)!} \frac{d^n}{dz^n} (z^{n-1} \ln(\xi(z)))|_{z=1} \qquad (1).$$

*for any non-negative integer n.*

Here $\xi(z)$ is the Riemann xi-function related with the Riemann $\zeta$-function by the well-known relation [1]

$$\xi(z) = \frac{1}{2} z(z-1) \pi^{-z/2} \Gamma(z/2) \varsigma(z) \qquad (2).$$

Two years later, Bombieri and Lagarias generalized Li's criterion [3]. Argumentation of [2, 3] goes roughly as follows. If $\rho = 1/2 + iT$, $T$ real and $i = \sqrt{-1}$, than $\left|\frac{\rho-1}{\rho}\right| = 1$ and hence can be written as $\exp(i\vartheta_i)$ where $\vartheta_i = \arctan \frac{T}{T^2 - 1/4}$. Let us introduce the sum $k_n = \Sigma_\rho (1 - (1 - \frac{1}{\rho})^n) = \Sigma_\rho (1 - \left(\frac{\rho-1}{\rho}\right)^n)$ over non-trivial Riemann function zeroes (*n* is non-negative integer, zeroes are counting taking into account their multiplicity, for *n*=1 contributions of complex conjugate zeroes should be paired when summing). For two complex conjugate "correct" Riemann function zeroes $\rho = 1/2 \pm iT$ we easily see that their contribution to sum $k_n$ is $2(1-\cos(n\vartheta_i))$, and hence non-negative; correspondingly, the sum $k_n$ is also non-negative. Quite the contrary, if some non-trivial Riemann function zero with $\text{Re}\,\rho \neq 1/2$ exists, for large enough *n* we will have an arbitrary large (by module) negative contributions from these zeroes, and it is straightforward to show that for infinitely many *n* this contribution can not be compensated by all other "correct" $1-\cos(n\vartheta_i)$ terms of the sum [3], whence infinitely many sums $k_n$ are to be negative.

This consideration immediately shows that the non-negativity of the sums $k_n$ is equivalent to the Riemann hypothesis. Li also demonstrated that these sums are equal to derivatives presented in eq. (1) (certainly, this is the most technically difficult part of his work; another derivation of this relation will be given shortly below).

## 2. GENERALIZED LI'S AND BOMBIERI-LAGARIAS CRITERIA

Now we note that for $\rho = 1/2 + iT$ and *any* real *a*
$\left|\dfrac{\rho - a}{\rho + a - 1}\right| = \left|\dfrac{-a + 1/2 + iT}{a - 1/2 + iT}\right| = 1$ and introduce the sum
$k_{n,a} = \Sigma_\rho (1 - \left(\dfrac{\rho - a}{\rho + a - 1}\right)^n) = \sum_\rho (1 - \left(1 - \dfrac{2a-1}{\rho + a - 1}\right)^n)$. To demonstrate that on RH all these sums are non-negative, just replace $\vartheta_i = \arctan\dfrac{T}{T^2 - 1/4}$ given above by $\vartheta_i = \arctan\dfrac{T(2a-1)}{T^2 - a^2 + a - 1/4}$ and repeat all the abovesaid. To demonstrate the inverse implication, let us briefly reproduce a slightly modified argument of Bomberi and Lagarias [3]; see their original paper for some more details.

Let $a<1/2$. We observe that for any Riemann zero $\rho = \sigma + iT$,
$\left|\frac{\rho-a}{\rho+a-1}\right|^2 = 1 + \frac{(1-2a)(2\sigma-1)}{|\rho+a-1|^2}$, and thus if $\sigma > 1/2$ we may find at least one zero for which $\left|\frac{\rho-a}{\rho+a-1}\right| > 1$. Because $\frac{(1-2a)(2\sigma-1)}{|\rho+a-1|^2}$ tends to zero when $|\rho_k|$ tends to infinity, maximum of this expression over $\rho$ is achieved and there are only finitely many, say $K$, zeroes $\rho_k$ for which $\left|\frac{\rho-a}{\rho+a-1}\right| = 1+t = \max$, for all others $\left|\frac{\rho-a}{\rho+a-1}\right| \leq 1+t-\delta$ for some fixed positive $\delta$. Clearly, taking $n$ large enough, the term $1 - \left(\frac{\rho_k-a}{\rho_k+a-1}\right)^n = 1-(1+t)^n \exp(in\vartheta_k)$ ($\vartheta_k$ is an argument of $\left(\frac{\rho_k-a}{\rho_k+a-1}\right)$) can be made very large by module and negative. Then, due to the Dirichlet's theorem on simultaneous Diophantine approximation, the sum of $1 - \left(\frac{\rho_k-a}{\rho_k+a-1}\right)^n$ over all $\rho_k$ can be made arbitrary close to $K(1-(1+t)^n)$ while the sum over all other zeroes is of the order of $O(n^2(1+t-\delta)^n)$, just due to their known density. The case $a>1/2$ is quite similar, so we have proven

Theorem 1. *Riemann hypothesis is equivalent to the non-negativity of sums* $k_{n,a} = \Sigma_\rho (1 - \left(\frac{\rho-a}{\rho+a-1}\right)^n) = \sum_\rho (1 - \left(1 - \frac{2a-1}{\rho+a-1}\right)^n)$ *taken over the Riemann xi-function zeroes for any real a, except a=1/2. Here n is non-negative integer, zeroes are counting taking into account their multiplicity, for n=1 contributions of complex conjugate zeroes should be paired when summing.*

Indeed, we proved this statement not only for the Riemann $\xi$-function zeroes but for certain multisets of complex numbers, see [3]. For completeness, here we formulate this result as a following

Theorem 2. (Generalized Bombieri – Lagarias' theorem). *Let a and $\sigma$ are arbitrary real numbers, $a < \sigma$, and R be a multiset of complex numbers $\rho$ such that*

(i) $\quad 2\sigma - a \notin R$

(ii) $\quad \sum_{\rho} (1+|\operatorname{Re}\rho|)/(1+|\rho+a-2\sigma|^2) < +\infty$

*Then the following conditions are equivalent*

*(a)* $\operatorname{Re}\rho \leq \sigma$ *for every* $\rho$;

*(b)* $\sum_{\rho} \operatorname{Re}(1 - \left(\dfrac{\rho - a}{\rho - 2\sigma + a}\right)^n) \geq 0$ *for n=1, 2, 3…*

*(c) For every fixed $\varepsilon > 0$ there is a positive constant $c(\varepsilon)$ such that*

$$\sum_{\rho} \operatorname{Re}(1 - \left(\dfrac{\rho - a}{\rho - 2\sigma + a}\right)^n) \geq -c(\varepsilon)e^{\varepsilon n}, \; n=1, 2, 3…$$

*If at the same conditions $a > \sigma$ is taken, the point (a) is to be changed to*

*(a')* $\operatorname{Re}\rho \geq \sigma$ *for every* $\rho$,

*points (b), (c) remain unchanged.*

Apparently the criterion *(ii)* can be recast in a more general form $\sum_{\rho} 1/(1+|\rho|^2) < +\infty$ [3], but this circumstance does not seem very important for the present author. The same point is relevant for the Theorem 3 (generalized Li's criterion) below.

As you see, the statement of Theorem 2 is formulated for *any* $\sigma$, not only for $\sigma = 1/2$, provided $\sigma \neq a$. To demonstrate this, just note that for $\rho = \sigma + iT$, $\left|\dfrac{\rho - a}{\rho + a - 2\sigma}\right| = \left|\dfrac{\sigma - a + iT}{a - \sigma + iT}\right| = 1$ and for $\rho = q + iT$,

$$\left|\dfrac{\rho - a}{\rho + a - 2\sigma}\right|^2 = 1 + \dfrac{4(\sigma - a)(q - \sigma)}{|\rho + a - 2\sigma|^2},$$ and then repeat all the abovesaid.

If, additionally to the aforementioned conditions of the generalized Bombieri – Lagarais' theorem, also the following takes place:

(iii) If $\rho \in R$, than $\bar{\rho} \in R$ with the same multiplicity as $\rho$

one can omit the operation of taking the real part in (b), (c), the expressions at question are real. (Here, as usual, $\bar{\rho}$ means a complex conjugate of $\rho$).

Following again the paper of Bombieri and Lagarais [3], we conclude this Section with the following

**Theorem 3. (Generalized Li's criterion).** *Let $a$ is an arbitrary real number, $a \neq \sigma$, and $R$ be a multiset of complex numbers $\rho$ such that*

(i)   $2\sigma - a \notin R$, $a \notin R$

(ii)   $\sum_{\rho}(1 + |\operatorname{Re}\rho|)/(1 + |\rho + a - 2\sigma|^2) < +\infty$, $\sum_{\rho}(1 + |\operatorname{Re}\rho|)/(1 + |\rho - a|^2) < +\infty$

(iii)   *If* $\rho \in R$, *than* $2\sigma - \rho \in R$

*Then the following conditions are equivalent*

(a) $\operatorname{Re}\rho = \sigma$ *for every* $\rho$;

(b) $\sum_{\rho} \operatorname{Re}\left(1 - \left(\dfrac{\rho - a}{\rho + a - 2\sigma}\right)^n\right) \geq 0$ *for any $a$ and $n=1, 2, 3...$*

(c) *For every fixed $\varepsilon > 0$ and any $a$ there is a positive constant $c(\varepsilon, a)$ such that* $\sum_{\rho} \operatorname{Re}\left(1 - \left(\dfrac{\rho - a}{\rho + a - 2\sigma}\right)^n\right) \geq -c(\varepsilon, a)e^{\varepsilon n}$, *for $n=1, 2, 3...$*

For clearly, in the conditions of the theorem we for all $\rho$ have $\mathrm{Re}\,\rho \leq \sigma$ and $\mathrm{Re}(2\sigma - \rho) \leq \sigma$, whence $\mathrm{Re}\,\rho = \sigma$. If, additionally to the aforementioned conditions of the generalized Li's criterion, also the following takes place:

(iv)  If $\rho \in R$, than $\bar{\rho} \in R$ with the same multiplicity as $\rho$

one can omit the operation of taking the real part in (b), (c), the expressions at question are real.

*Remark.* Similarly to the Li's criterion, generalized Li's criterion can be applied also to numerous other zeta-functions, as this was shown first by Li himself for Dedekind zeta-function [2], and afterwards was the subject of a number of sequel papers by other authors. We will not pursue this line of researches here. Similarly, at a moment we put aside the questions concerning the relation of these generalized criteria with so-called Weil's explicit formula for the theory of prime numbers [4], see [3] and [5].

Our next aim is to establish relation "of the Li's type" similar to eq. (1), viz. the relation between sums $k_{n,a}$ and certain derivatives of the Riemann xi-function. For this we will use the generalized Littlewood theorem concerning contour integrals of logarithm of an analytical function, recently used by us to establish numerous equalities equivalent to the Riemann hypothesis [6], which we reproduce below. The proof [6] is a straightforward modification of familiar and well known corresponding Littlewood theorem (or Lemma) proof, see e.g. [7]. Actually, this theorem has been more or less explicitly used in Riemann researches already by Wang who in 1946 established the first integral equality equivalent to the Riemann hypothesis [8].

**Theorem 4 (Generalized Littlewood theorem).** *Let C denotes the rectangle bounded by the lines $x = X_1$, $x = X_2$, $y = Y_1$, $y = Y_2$ where $X_1 < X_2$, $Y_1 < Y_2$ and let f(z) be analytic and non-zero on C and meromorphic inside it, let also g(z) is analytic on C and meromorphic inside it. Let F(z)=ln(f(z)), the logarithm being defined as follows: we start with a particular determination on $x = X_2$, and obtain the value at other points by continuous variation along y=const from $\ln(X_2 + iy)$. If, however, this path would cross a zero or pole of f(z), we take F(z) to be $F(z \pm i0)$ according as we approach the path from above or below. Let also the poles and zeroes of the functions f(z), g(z) do not coincide.*

$$\text{Then } \int_C F(z)g(z)dz = 2\pi i(\sum_{\rho_g} res(g(\rho_g)\cdot F(\rho_g))) - \sum_{\rho_f^0} \int_{X_1+iY_\rho^0}^{X_\rho^0+iY_\rho^0} g(z)dz + \sum_{\rho_f^{pol}} \int_{X_1+iY_\rho^{pol}}^{X_\rho^{pol}+iY_\rho^{pol}} g(z)dz$$

*where the sum is over all $\rho_g$ which are poles of the function g(z) lying inside C, all $\rho_f^0 = X_\rho^0 + iY_\rho^0$ which are zeroes of the function f(z) counted taking into account their multiplicities (that is the corresponding term is multiplied by m for a zero of the order m) and which lye inside C, and all $\rho_f^{pol} = X_\rho^{pol} + iY_\rho^{pol}$ which are poles of the function f(z) counted taking into account their multiplicities and which lye inside C. For this is true all relevant integrals in the right hand side of the equality should exist.*

*Remark.* Actually, the case of the coincidence of poles and zeroes of the functions *f(z), g(z)* often does not pose real problems and can be easily considered. We have dealt with a few such cases before [6].

The subtle moment related with this generalized Littlewood theorem is the circumstance that the function *arg(F(z))* (imaginary part of the *ln(f(z))*) is *not* continuous on the left border of the contour (segment $X_1+iY_1$, $X_1+iY_2$) if there are zeroes or poles of the function *f(z)* inside the contour. This is explicitly stated in the theorem condition: *If, however, this path would cross a zero or pole of f(z), we take F(z) to be $F(z \pm i0)$ according as we approach the path from above or below.* In practice, this means that when calculating the corresponding part of the contour integral, viz. the integral $-\int_{X_1+iY_1}^{X_1+iY_2} \arg(F(z))g(z)dz$ (minus sign comes from the necessity to round the contour counterclockwise), $\pm 2\pi i l$ jumps should be added to an argument function at a point $X_1+iY_{z,p}$ whenever a zero or a pole of an order *l* of the function *f(z)* occurs somewhere at a point $X+iY_{z,p}$ inside the contour. Corresponding integral should be properly modified if the use of a continuous argument branch is desirable. See our paper [9] for details, we also would like to note that the appropriateness of the necessary modification of an argument has been numerically tested (and confirmed) by us for a number of integrals, e.g. for the integral $\int_0^\infty \frac{t \arg(\varsigma(1/4+it))}{(1/16+t^2)^2}dt = \pi\frac{\varsigma'(1/2)}{\varsigma(1/2)} - 9\pi - \pi \sum_{\rho,\sigma_k > 1/4, t_k > 0}\left(\frac{1}{t_k^2+1/4}\right)$ (similar equality in the form $\int_0^\infty \frac{t \arg(\varsigma(1/2+it))}{(1/16+t^2)^2}dt = \pi\frac{\varsigma'(3/4)}{\varsigma(3/4)} - \frac{32\pi}{3}$ is equivalent to the Riemann hypothesis, our Theorem 5 from [6]). However, for what follows the asymptotic of the function *g(z)* for large values of $X_1$ tending to minus infinity makes this modification irrelevant, the value of the integral $-\int_{X_1+iY_1}^{X_1+iY_2} \arg(F(z))g(z)dz$ tends to zero anyway.

First, as an exercise, we use this theorem to establish the Li's relation (1). For this, let us consider the rectangular contour $C$ with vertices at $\pm X \pm iX$ with real $X \to +\infty$, if some Riemann zero occurs on the contour just shift it a bit to avoid this, and consider a contour integral $\int_C g(z)\ln(\xi(z))dz$ where

$$g(z) = \frac{n}{(z-1)^2}\left(\frac{z}{z-1}\right)^{n-1} - \frac{n}{(z-1)^2} \qquad (3).$$

Known asymptotic of the logarithm of the xi-function for large |z|, $\cong O(z\ln z)$ guaranties the "disappearance" of the contour integral value (it tends to zero when $X \to \infty$ due to the asymptotic $g(z) \cong O(1/z^3)$) thus we get, after division by $2\pi i$,

$$n\frac{1}{(n-1)!}\frac{d^n}{dz^n}(z^{n-1}\ln(\xi(z)))|_{z=1} - n\frac{\xi'}{\xi}(1) - \sum_\rho (1-(1-\frac{1}{1-\rho})^n) - n\sum_\rho \frac{1}{\rho} = 0 \quad (4).$$

(Complex conjugate zeroes are to be paired when calculating $\sum_\rho \frac{1}{\rho}$ and $\sum_\rho (1-(1-\frac{1}{1-\rho})^n)$ for n=1). Here the first term is the contribution of the n+1 order pole of g(z) at z=1, second term is the contribution of the second order pole arising from the term $-\frac{n}{(z-1)^2}$ occurring in (3), the third and fourth terms are the integrals $-\int_{-\infty+iT_i}^{\rho_i} g(z)dz$. (Clearly,

$\frac{n}{(z-1)^2}\left(\frac{z}{z-1}\right)^{n-1} = \frac{d}{dz}(1-\left(1-\frac{1}{1-z}\right)^n)$ which explains while function g(z) in form (3) is used; the term $-\frac{n}{(z-1)^2}$ is added just to ensure the asymptotic $g(z) \cong O(1/z^3)$ necessary to bring the contour integral value to zero. Note also evident $\sum_\rho (1-(1-\frac{1}{1-\rho})^n) = \sum_\rho (1-(1-\frac{1}{\rho})^n)$ ). We know that $\frac{\xi'}{\xi}(1) = -\sum_\rho \frac{1}{\rho}$ [1], so that we have $n\frac{1}{(n-1)!}\frac{d^n}{dz^n}(z^{n-1}\ln(\xi(z)))|_{z=1} = \sum_\rho (1-(1-\frac{1}{\rho})^n)$ which is a relation we have searched for.

Quite similar consideration is applied to our case, where now we introduce the function

$$\tilde{g}(z) = -\frac{n(2a-1)(z-a)^{n-1}}{(z+a-1)^{n+1}} + \frac{n(2a-1)}{(z+a-1)^2} \qquad (5)$$

and consider contour integral $\int_C \tilde{g}(z)\ln(\xi(z))dz$ taken round the same contour as above. Application of Theorem 4 (generalized Littlewood theorem) gives

$$-\frac{n(2a-1)}{(n-1)!}\frac{d^n}{dz^n}((z-a)^{n-1}\ln(\xi(z)))|_{z=1-a} + n(2a-1)\frac{\xi'}{\xi}(z)|_{z=1-a} - \sum_\rho (1-(\frac{\rho-a}{\rho+a-1})^n) + n(2a-1)\sum_\rho \frac{1}{\rho+a-1} = 0 \qquad (6)$$

(Again, complex conjugate zeroes are to be paired whenever necessary). Using well known $\frac{\xi'}{\xi}(z)|_{z=1-a} = -\sum_\rho \frac{1}{\rho+a-1}$ [1] and reminding our theorem 1 we have

Theorem 5. *Riemann hypothesis is equivalent to the non-negativity of all derivatives* $\frac{1}{(n-1)!}\frac{d^n}{dz^n}((z-a)^{n-1}\ln(\xi(z)))|_{z=1-a}$ *for all non-negative integers n and any real a<1/2; correspondingly, it is equivalent also to the non-positivity of all derivatives* $\frac{1}{(n-1)!}\frac{d^n}{dz^n}((z-a)^{n-1}\ln(\xi(z)))|_{z=1-a}$ *for all non-negative integers n and any real a>1/2.*

*Remark.* Another possibility to arrive to the same conclusions is to see the formula $\left|\dfrac{\rho-a}{\rho+a-1}\right|=\left|\dfrac{-a+1/2+iT}{a-1/2+iT}\right|=1$ as a precursor for the conformal mapping $s=\dfrac{z-a}{z+a-1}$. For $a<1/2$ and $\operatorname{Re}z\leq 1/2$, module of $s$ is always less or equal to 1; this equality is realized only on the line $z=1/2+it$. Correspondingly, on RH the function $\ln\xi(\dfrac{z-a}{z+a-1})$ is analytic in the interior of the discus $|s|<1$. We will not pursue this line of researches here, see again [2, 3] and our paper [8] where similar idea was used to generalize Balazard – Saias – Yor's criterion equivalent to the Riemann hypothesis [9]. Similarly, for more general case $s=\dfrac{z-a}{z+a-2\sigma}$, if $a<\sigma$ and $\operatorname{Re}z\leq\sigma$, module of $s$ is always less or equal to 1; and if $a>\sigma$ and $\operatorname{Re}z\geq\sigma$, this module also is always less or equal to 1. This illustrates again our Theorem 2 (the generalized Bombieri – Lagarias' theorem).

*Remark.* Along similar lines, analogous formulae connecting generalized Li's sums and certain derivatives of the logarithm, can be established for numerous other zeta-functions. We will not pursue this line of researches here.